\documentclass[11pt]{amsart}
\numberwithin{equation}{section}
\newtheorem{theorem}{Theorem}
\newtheorem{definition}[theorem]{Definition}
\newtheorem{proposition}[theorem]{Proposition}

\newtheorem{corollary}[theorem]{Corollary}
\numberwithin{theorem}{section}
\newtheorem{remark}[theorem]{Remark}
\def\R{\bf R}

\def\al{\aligned}
\def\eal{\endaligned}
\def\M{{\bf M}}
\def\be{\begin{equation}}
\def\ee{\end{equation}}
\def\lab{\label}
\def\a{\alpha}

\def\e{\epsilon}

\def\R{{\bf R}}
\def\lam{\lambda}

\def\al{\aligned}
\def\p{\partial}
\numberwithin{equation}{section}

\usepackage[backref=page]{hyperref}

\begin{document}

\tracingpages 1
\title[rigidity]{\bf A rigidity result for ancient Ricci flows}
\author{Qi S. Zhang}
\address{Department of
Mathematics  University of California, Riverside, CA 92521, USA }
\date{June 2024,  MSC2020: 58J35.}

\begin{abstract}
Using a size condition of the sharp  log Sobolev functional (log entropy) near infinity only, we prove a rigidity result for ancient Ricci flows without sign condition on the curvatures. The result is also related to the problem of identifying type II ancient Ricci flows and their backward limits.

\end{abstract}
\maketitle
\section{Statement of result}

Let $(M, g(t))$, $t \in (-\infty, 0)$ be an ancient solution to the Ricci flow on a noncompact Riemannian manifold $M$ of dimensional $n \ge 3$, which is $\kappa$ non-collapsed and is of bounded curvature. In this paper, we prove that if a sharp version of the $W$ entropy at infinity alone is close to $0$, then $(M, g(t))$ is flat. This includes the case when $(M, g(t))$  is $C^2$ asymptotically close to a round cone near infinity uniformly in $t$, whose cone angle is sufficiently close to $2 \pi$, i.e. that of the Euclidean space.
In particular, uniformly asymptotically flat manifolds are ruled out as type II singularity models of the Ricci flow.
The restriction on the entropy at infinity or the cone angle at infinity is reasonable and necessary in general, at least qualitatively.
This assertion is valid in view of Gap Lemma 3.1 in the paper \cite{An:1} by Anderson on flatness of Ricci flat manifolds with almost Euclidean volume and examples of non-flat steady gradient solitons such as the Eguchi-Hansen metric. Recall that above named metric is Ricci flat but not flat and the volume of large geodesic balls are in the same ordered as the Euclidean ones but not asymptotically close. In fact, Corollary \ref{corig} (a)  below can be regarded as an extension of the above result in \cite{An:1} to gradient Ricci soliton cases. The main difference is that no vanishing condition on the Ricci curvature is required at the expense of some assumptions on the metric near infinity.

 Another motivation for this kind of results come from  the study of Ricci flows, for which identifying non-flat ancient solutions, including gradient Ricci solitons, is a central problem since many of them  form blow up limits of singularities, i.e. singularity models.
 In dimension 3, Perelman \cite{P:1} proved that a backward in time limit of  ancient $\kappa$ solutions are Ricci solitons. Recently, Brendle \cite{Br:1}  further proved that type II noncompact ancient $\kappa$ solutions are the Bryant solitons. This implies, with previous and further work, that ancient $\kappa$ solutions are either the Bryant solitons or shrinking cylinders or their quotients.  The proofs rely heavily on the property that the sectional curvatures are nonnegative.

 In higher dimensions, a similar convergence to non-flat solitons for a type I $\kappa$ non-collapsed ancient solution or singularity models was proven some years ago in \cite{EMT:1} and \cite{CZ:1} independently. Recall that an ancient solution $(M, g(t))$ is type I if the curvature tensor satisfies $|Rm_{g(t)}(\cdot, t)| \le C/|t|$.
  Recently, for a finite time singularity model which is partially type I,  we proved that a blow up limit is a
gradient shrinking soliton \cite{Z22}.
  In general a similar statement on backward limits, when shrinking solitons is replaced by steady solitons,  expected by workers in the field,  is still lacking for a type II (non type I) ancient solution or singularity models when the dimension is 4 or higher. A helpful step in solving the problem is to identity examples of type II ancient solutions or singularity models for the Ricci flow. One recent example can be found in \cite{Ap}, where it was shown that the Eguchi-Hanson manifold is such an example.  Theorem \ref{thmrig} and Corollaries below can also be seen as a small step in this investigation by ruling out certain type II ancient solutions from manifolds mentioned in the first paragraph.

Before going into details, let's introduce notations, concepts and definitions to be used in the paper.
We use $M$ to denote a $n (\ge 3)$ dimensional complete  Riemannian manifold and $g(t)$ to denote the
metric at time $t$; $d(x, y, t)$ is the geodesic distance under $g(t)$; Unless stated otherwise, we
assume the curvature tensor is bounded at each time $t$.
$B(x, r, g(t)) = \{ y \in {\M} \ | \ d(x, y, t) < r \}$  is the geodesic ball of radius $r$, under metric $g(t)$,
centered at $x$; when no confusion arises we may also use $B(x, r)$ or $B(x, r, t)$
to denote $B(x, r, g(t)) $;  and $|B(x, r, t)|_{g(t)}$ is the volume of
$B(x, r, t)$ under $g(t)$;    $dg(t)$ is the volume element; $x_0$ or $0$ is a reference point on $M$. We also
reserve $R=R(x, t)$ as the scalar curvature under $g(t)$. $\Delta$ and $\nabla$ are the Laplace-Beltrami operator and gradient with respect to a metric and if there is a need, the corresponding metric is denoted by a subscript such as $\Delta_g$ e.g..  A generic positive constant is denoted
by $C$ or $c$ whose value may change from line to line. When we say that a sequence of
pointed manifolds converges in $C^\infty_{loc}$ sense, we mean they converge in the usual
Cheeger-Gromov sense. That is, subject to diffeomorphisms, the metrics converge
in $C^\infty_{loc}$ sense.

Let us recall the scaling invariant sharp log Sobolev functionals, originally  introduced by Weissler \cite{W:1} on $\R^n$, which can be found in \cite{Z14} e.g.

\begin{definition} ( Sharp Log Sobolev functional, infimum, infimum at
infinity)
 \lab{deflv} Let $(M, g)$ be a $n$ dimensional Riemannian
manifold with metric $g$ and $D \subset M$ be a domain. Suppose the $F$ functional is nonnegative, i.e.
\[
\inf_{\Vert v \Vert_{L^2}=1}  \underbrace{\int_D ( 4 |\nabla v |^2 + R v^2 ) dg}_{F(v)} \ge 0.
\]

(a). Given  functions $ v \in W^{1, 2}_0(D, g)$ with $\Vert v
\Vert_{L^2(D)}=1$,   the sharp log Sobolev
functionals is defined by
\be
\lab{lvg}
\al
 L(D, v, g) &= - \underbrace{\int_D v^2 \ln v^2 dg}_{N(v)} + \frac{n}{2} \ln
\left(\int_D ( 4 |\nabla v |^2 + R v^2 ) dg \right) + s_n \\
&\equiv - N(v) +   \frac{n}{2} \ln F(v) + s_n.
\eal
\ee Here $s_n=-(n/2) \ln (2 \pi e n)$.

(b). The infimum of the sharp log Sobolev functional is denoted by
\[
\lam(D, g) =
\inf \{ L(D, v, g) \, | \,  v \in W^{1, 2}_0(D, g), \quad \Vert v \Vert_{L^2(D)}=1 \}.
\]

(c). The infimum of the sharp log Sobolev
functional at infinity is
\[
\lam_\infty(M, g) = \lim_{\rho \to \infty} \lam(M-B(x_0, \rho), g)
\]where $x_0$ is a reference point in $M$.

\end{definition}

\begin{remark}  The sharp log Sobolev functional is a
scaling invariant version of the log Sobolev functional with parameters originally
introduced by Gross \cite{G:1} and Federbush \cite{F:1}. It can also be regarded as the infimum of Perelman's W entropy over all scales \cite{P:1}. Part (c) is fashioned from an idea in \cite{Lio:1} by P.L. Lions.

When $F(v)$ becomes $0$ but $N(v)$ is finite, the functional $L$ is
regarded as $-\infty$.
\end{remark}

\begin{definition} (gradient Ricci solitons)  A Riemannian manifold $(M, g)$ is called a
gradient Ricci soliton if there exists a smooth function $f$ on
$M$ and a constant $\e_*$ such that
\be
\lab{defGRSol}
Ric+ Hess f +
\frac{\e_*}{2} g =0.
 \ee

$(M, g)$ is called a expanding, steady and shrinking gradient Ricci soliton if
$\e_*>0, \e_*=0$ and $\e_*<0$ respectively.
\end{definition}

\begin{definition} Let $\kappa$ be a positive number.
A Riemannian manifold $(M, g)$ of dimension $n$ is called $\kappa$ non-collapsed at scale $r_*$ if for any $r \in (0, r_*]$ and $x \in M$ the following holds.  If $|Rm| \le 1/r^2$ in $B(x, r)$, then $|B(x, r)| \ge \kappa r^n$.
\end{definition}

\begin{definition}
\lab{defAF} Given a number $\a \in (0, 1]$ and $n \ge 3$,
let $(\R^n, g_\a)$ be the standard round cone whose metric is given by $g_\alpha = d^2 r + \alpha^2 r^2  g_{S^{n-1}}$, where $ g_{S^{n-1}}$ is the standard metric in the unit $n-1$ sphere $S^{n-1}$.
A complete, noncompact Riemannian manifold $M$ is called  {\bf Asymptotically $\e$ close to a round $\a$-cone} if
there is a partition $M=M_0 \cup M_\infty$, which satisfies the following properties.

(i). $M_0$ is compact with a reference point $0$.

(ii). $M_\infty$ is the disjoint union of finitely many components each of which is diffeomorphic to $({\R^n} - B(0, r_0))$ for some $r_0 \ge 1$.

(iii). Under the coordinates induced by the diffeomorphism,  the metric $g_{ij}$ satisfies, for $x \in
M_\infty$, $\e>0$
\[
\al
&|g_{ij}(x) -  (g_\a)_{ij}(x)| \le \e, \quad |\partial_k
g_{ij}(x)- \partial_k
(g_\a)_{ij}(x)| \le  \e |x|^{-1}, \\
& |\partial_k  \partial_l
g_{ij}(x)-\partial_k  \partial_l
(g_\a)_{ij}(x)| \le  \e |x|^{-2}.
\eal
\]
\end{definition}

\begin{remark}
\lab{re2.1}
For convenience we will assume  that $M_\infty$ has only one connected
component. This assumption does not reduce any generality.
\end{remark}

When $\a=1$, the above definition contains, as special cases,  the usual Asymptotically Flat (AF) manifolds for which the order of decay for the gradient and second derivative of $g_{ij}$ is faster than $-1$ and $-2$ respectively. Due to Theorem (1.1) in \cite{BKN:1}, if $M$ has one end,  the curvature tensor decays sufficiently fast near infinity and
$|B(0, r)| \ge c r^n$ when $r$ is large, then $M$ is AF manifold. Here $n$ is the dimension.

\medskip

We are now ready to state the main results of the paper.

\begin{theorem}
\lab{thmrig}

There exists a positive number $\e_0$, depending only on the dimension $n (\ge 3)$ such that the following is true. Let $(M, g(t), x_0)$,  $\partial_t g_{ij} = - 2 R_{ij}$, $t \in (-\infty, 0]$ be a complete, noncompact  Ricci flow
 with bounded curvature,  which is $\kappa$ non-collapsed at scale $1$. If, for some $r_0 \ge 1$, the inequality
 \be
 \lab{lambbce0}
 \lam(M-B(x_0, r_0, g(t)), g(t)) \ge -\e_0
 \ee holds for all $t \le 0$, then $(M, g(t))$ is the standard $\mathbf{R^n}$. The converse is also true.
\end{theorem}

\begin{corollary}
\lab{corig}
There exists a positive constant $\alpha_0 \in (0, 1)$, depending only on the dimension $n (\ge 3)$, and another sufficiently small constant $\e$ depending only on $n, \a_0$ such that the following conclusions hold.

(a).  Suppose $(M, g)$ is asymptotically $\e$ close to a round $\a$-cone with $\a \ge \a_0$ and $R \ge 0$. Then $(M, g)$ is not a non-flat, gradient steady or expanding soliton.

(b).
Let $(M, g(t))$,  $\partial_t g_{ij} = - 2 R_{ij}$, $t \in (-\infty, 0]$ be a complete, noncompact  Ricci flow
 with bounded curvature, which is $\kappa$ non-collapsed at scale $1$. Suppose $(M, g(t), 0)$ is asymptotically $\e$ close to a round $\a$-cone with $\a \ge \a_0$ uniformly in time. i.e. the radius $r_0$ in Definition \ref{defAF} is independent of time.
 Then $(M, g(t))$ is flat.

\end{corollary}

\begin{remark}

(a). It should be noted that the asymptotic conditions are imposed near infinity instead of on the whole manifold. In fact the conclusion of the theorem will follow relatively quickly if one imposes a similar condition on the whole manifold. See \cite{Yo} by Yokota e.g., and Proposition \ref{prAFsob} below, which is the starting point for the proof of the theorem.   There are many examples of non-flat shrinking, expanding  gradient solitons on asymptotically conic manifolds whose cone angle is far from Euclidean. See \cite{FIK:1} and \cite{Ca:1} e.g. Hence the restriction on the parameter $\a$ in the theorem is necessary qualitatively.

(b).  Since complete noncompact flat manifolds with maximum volume growth are Euclidean spaces with the standard metric, the corollary implies that ancient solutions satisfying the given conditions are Euclidean.

(c). In particular the conclusions on the theorem hold when $(M, g)$ is the standard asymptotically flat manifold for which  $\a=1$ and which has faster decay for the covariant derivatives of the metric. This class is interesting due to connections
to general relativity. Useful properties of these kind of Ricci
flows have bee proven in \cite{DM:1}, \cite{OW:1}. For example,
they proved that the AF property is preserved under Ricci flow. In the AF case, Part (a) of the Corollary seems to extend Theorem A.3 in \cite{Liy:1} by Yu Li, where the decay conditions is some order faster than condition (iii) of Definition \ref{defAF}.
 This extra decay is crucial for many results on the subject but is not needed here.
See also \cite{Z18} for related results in the stationary metric case.

(d).  If one assumes the Ricci curvature is nonnegative, then it is well known that non-flat gradient shrinking solitons must have $0$ asymptotic volume ratio, namely $|B(0, r)|/r^n \to 0$ as $r \to \infty$. See \cite{CN:1} e.g. Also part (a) of the theorem does not hold when the dimension is 2. See the 2 dimensional expanding soliton example in Sec.5, Chapter 4 of \cite{CLN:1}.
\end{remark}

As a by product, we also obtained the following necessary and sufficient condition on existence of minimizers for the sharp log Sobolev functional on solitons.

\begin{proposition}
\lab{prminiSLS}
Let $(M, g)$ be a complete, noncompact Ricci soliton of dimension $n \ge 3$ which is $\kappa$ non-collapsed at scale one.
Suppose the curvature tensor is bounded. Then the sharp log Sobolev functional has a minimizer in $W^{1, 2}(M)$ if and only if $(M, g)$ is a gradient shrinking soliton.
\end{proposition}

Let us outline the proof of the theorem. The starting point is a rigidity result involving the infimum of the sharp log Sobolev functional $\lam(M, g)$ , which states that $\lam(M, g)=0$ if and only if $M=\R^n$ with the usual metric. This will imply that for  non-flat ancient solutions, the corresponding $\lam(M, g(t))$ will be bounded from above by a negative number $-\delta_*$ as $t \to -\infty$. Under the asymptotic assumption in the theorem, $\lam_\infty(M, g(t))$ will be strictly greater than $\lam(M, g(t))$ when $t \to -\infty$. According to the main result in \cite{Z14}, using P. L. Lions' concentrated compactness method \cite{Lio:1}, a minimizer for $\lam(M, g(t))$ can be found. See also
 \cite{DE:1}  where Dolbeault and Esteban treated a similar
functional on the cylinder $S^n \times \R$.
 Once a minimizer is found at a time level $t$, we still need to prove that it will not fizzle to $0$ as $t \to -\infty$. This is the main technical work of the paper, which is based on the compensated compactness method again. Now we can use Perelman's
monotonicity formula to show that $(M, g(t))$ is a shrinking gradient Ricci
soliton as $t \to -\infty$ since $\lim_{t \to -\infty} \lam(g(t))$ exists. We then will prove, with the help of a classical result in \cite{T:1},  that if a manifold is both a shrinking and steady or expanding soliton, then it must be flat.

\section{Proof of the theorem}

Our starting point is the following result on the sharpness and rigidity of the best constant of the sharp log entropy $L$. The first statement (a) is in a similar spirit to the main result in \cite{Le:1} by M. Ledeux, assuming the Ricci curvature is non-negative.  The difference is that we do not need any curvature assumption except its boundedness.
The second statement (b) is similar to the main result \cite{Yo} by Yokota where the relevant criteria is on Perelman's reduced volume.

\begin{proposition}
\lab{prAFsob} (a). Let $(M, g)$ be a complete Riemannian manifold of dimension $n \ge
3$ with bounded curvature.
Then
\[
\lambda(M, g)=0
\]if and only $M$ is isometric to $\R^n$ with the standard Euclidean metric. In other words, if the sharp log Sobolev inequality holds:
\[
 \int_D v^2 \ln v^2 dg \le \frac{n}{2} \ln
\left(\int_D ( 4 |\nabla v |^2 + R v^2 ) dg \right) + s_n, \quad \forall v \in W^{1, 2}(M), \quad \Vert v \Vert_{L^2}=1,
\]then $M$ is isometric to $\R^n$ with the standard Euclidean metric.

(b). Let $(M, g(t))$ be a nonflat ancient solution  of the Ricci flow with bounded curvature. Then, there exists a positive number $\delta_*$, depending only on the dimension $n$  such that
\[
\lim_{t \to - \infty} \lambda(M, g(t)) \le -\delta_*.
\]

(c). Let $(M, g)$ be a  time slice of a nonflat expanding soliton with bounded curvature which is also $\kappa$ non-collapsed at all scales. Then, there exists a positive number $\delta_*$, depending only on the dimension $n$  such that
\[
 \lambda(M, g) \le -\delta_*.
\]

\end{proposition}

\proof (a). One direction of the result is the well known fact that the best constant of the sharp log Sobolev inequality on $\R^n$ is $0$.

For the other direction, we assume $\lambda(M, g)=0$. This implies that the W entropy is bounded from below since $\lambda(M, g)$ is the infimum of the $W$ entropy for all scales , which shows that $M$ is $\kappa$ non-collapsed at scale one at least. In fact, it is known, as described several lines below, that $M$ is $\kappa$ non-collapsed at all scales.
By Proposition 2.4 in \cite{Z14}, the manifold $M$ is a shrinking gradient Ricci soliton.
Notice that we do not assume the scalar curvature is nonnegative to begin with. But this property follows from \cite{Ch:1} once we know $M$ is a shrinking soliton. Now we can conclude that $M$ is isometric to $\R^n$ using Corollary 1.1 (3) in \cite{Yo}, whose condition includes as a special case  the  condition that the W entropy at scale $1$ is greater than a small negative number. See Remark 6.5 in that paper and also a later paper \cite{Liy:1}. Since $\lambda(M, g)$ is the infimum of the $W$ entropy for all scales, our condition that $\lambda(M, g)=0$ implies that the $W$ entropy at scale $1$ is at least $0$ which is covered by the corollary cited above.

(b).
Suppose for contradiction that the conclusion is false. We can then find a sequence of non-flat, $\kappa$ non-collapsed  ancient solutions $(M_k, g_k(t))$ with bounded curvature such that, for some $t_k \le 0$,
\be
\lab{lamkdk0}
\lim_{t_k \to -\infty} \lam(M_k, g_k(t_k)) \equiv -\delta_k \to 0, \quad k \to \infty.
\ee  It is known that for fixed $k$, $\lam(M_k, g_k(t))$ is monotone non-decreasing in time $t$, which follows from the monotonicity of the $W$ entropy. So the above limit is well defined.

Note that $|Rm_{g_k}|$ may not be uniformly bounded for all $k$.
For each $k$, we choose a point $y_k \in M_k$ and time $s_k \le t_k$ such hat
\be
\lab{alk12}
\alpha_k \equiv | Rm_{g_k}(y_k, s_k)| \ge \frac{1}{2} \sup_{x \in M_k, t \le 0} | Rm_{g_k}(x, t)|.
\ee Consider the scaled ancient solutions
\[
\widetilde{g}_k = \alpha_k g_k(\cdot, \alpha^{-1}_k s+ s_k), \quad s \le 0.
\]By our choice of $\alpha_k$, the norm of the curvature tensors for $(M_k, \widetilde{g}_k(\cdot, s), y_k)$ are uniformly bounded by $2$ for all $k,  s \le 0$.

Due to the scaling invariance property of $\lam$ and \eqref{lamkdk0} and the monotonicity in time, we also have
\be
\lab{lamkwank}
\lam(M_k, \widetilde{g}_k(0)) = \lam(M_k, g_k(s_k)) \ge \lim_{t \to -\infty} \lam(M_k, g_k(t)) = -\delta_k \to 0.
\ee  Since $\delta_k \to 0$, we can assume, without loss of generality that $\delta_k \le 1$.  We argue that $(M_k, g_k(t))$
   is uniformly $\kappa$ non-collapsed at all scale for all space time point. i.e., for one constant $\kappa>0$,
\[
|B(x, r)| \ge \kappa r^n, \quad \text{if} \quad |Rm| \le 1/r^2 \quad \text{in} \quad B(x, r)
\]for all $x \in M_k$, $r>0$ and $t \le 0$. Here for simplicity we have suppressed the dependence of geometric quantities on the metrics $g_k(t)$. Indeed, by monotonicity of $\lam$ in $t$, we have
\[
\lam(M_k, g_k(t)) \ge -1,\qquad  t \le 0, \quad \forall k=1, 2, 3, ...
\]This implies that the $W$ entropy at any scale for $(M_k, g_k(t))$, which is not smaller than $\lam(M_k, g_k(t))$ is also bounded from below by $-1$. Then the $\kappa$ non-collapsing property follows by known procedure. See Proposition 2.6 in \cite{Z14} e.g.

By scaling invariance, the scaled flows $(M_k, \widetilde{g}_k(s))$ are also $\kappa$ non-collapsed at all scales uniformly for all $k$ and $s \le 0$. This and the boundedness of the curvatures imply, by Hamilton's compactness theorem for Ricci flows, that  a subsequence of the pointed flows $\{(M_k, \widetilde{g}_k(s), y_k)\}$, converges, in $C^\infty_{loc}$ topology,  to a pointed ancient solution $(M_\infty, \widetilde{g}_\infty(s), y_\infty)$.
By \eqref{alk12}, we deduce
\be
\lab{rm>1/2}
|Rm_{g_\infty}(y_\infty, 0)| \ge 1/2
\ee and by \eqref{lamkwank}, we infer
\[
\lam(M_\infty, g_\infty(0)) \ge \limsup_{k \to \infty} \lam(M_k, \widetilde{g}_k(0)) \ge
\medskip0.
\]We mention that the first inequality in the preceding line can be proven using a similar argument as Lemma 6.28 in \cite{chowetc3} where such inequality was established for the infimum of the $W$ entropy for all parameters.  Therefore $\lam(M_\infty, g_\infty(0))=0$ which, by part (a), implies that $(M_\infty, g_\infty(0))$ is isometric to $\R^n$ with the standard metric, reaching a contradiction with \eqref{rm>1/2}. This completes the proof of part (b).

(c). The key is the observation that for expanding solitons with bounded curvature, the $i$-TtH order covariant derivatives of the curvature tensor are bounded also. Indeed we can solve the Ricci flow $(M, g(t))$ with initial metric $g(0)=g$. Since $(M, g)$ is an expanding soliton, we know that $g(t)=(1+t) \phi^*_t g(0)$ where $\phi_t$ is a one parameter family of diffeomorphisms.
Let $A=\sup |Rm_{g(0)}|_{g(0)}$. Then
\[
\sup |Rm_{g(t)}|_{g(t)} \le A, \quad \forall t \ge 0.
\]According to Shi's gradient estimates \cite{Sh:1}, for any $t \in [1, 2]$, we have, for some positive constant $C_i$, that
\[
|\nabla^i Rm_{g(t)}(\cdot, t) |_{g(t)} \le C_i A, \quad  t \in [1, 2].
\]Scaling back to the initial metric $g(0)$, which is just the original $g$, we deduce
\be
\lab{dkrm}
|\nabla^i Rm_g |_{g} \le 2^{i} C_i A.
\ee  With this estimate in hand, the rest of the proof is similar to part (b).

Suppose for contradiction that the conclusion is false. We can then find a sequence of non-flat, $\kappa$ non-collapsed  expanding solitons $(M_k, g_k)$ with bounded curvature such that,
\be
\lab{lamkdk2}
\lim_{k \to \infty} \lam(M_k, g_k) \equiv -\delta_k \to 0, \quad k \to \infty.
\ee  Since $\lam(M_k, g_k)$ is scaling invariant, we can make a scaling if necessary to ensure that
\[
\sup |Rm_{g_k}|_{g_k} =1.
\]Due to the $\kappa$ non-collapsing assumption and \eqref{dkrm} with $g$ replaced by $g_k$ and $A=1$, we can again extract a subsequence,  denoted by $(M_k, g_k, x_k)$, which converges in $C^\infty_{loc}$ sense to a manifold $(M_\infty, g_\infty, x_\infty)$. Here $x_k$ is a point such that $|Rm_{g_k}(x_k)|_{g_k} \ge 1/2$.  Moreover $\lam(M_\infty, g_\infty)=0$ and hence it is the Euclidean space with the standard metric by part (a) of the proposition. But this contradicts with
\[
\sup |Rm_{g_\infty}|_{g_\infty} \ge \lim_{k \to \infty} |Rm_{g_k}(x_k)|_{g_k} \ge 1/2,
\] completing the proof of part (c).

\qed

The next result concerns size estimate of $\lam$ and $\lam_\infty$ of a round cone and manifolds which are asymptotically $\e$ close to round cone.

\begin{proposition}
\lab{prlamcaaca}
(a).
For a parameter $\alpha \in (0, 1]$, let $(\R^n, g_\a)$ be the standard round cone whose metric is given by $g_\alpha = d^2 r + \alpha^2 r^2  g_{S^{n-1}}$, where $ g_{S^{n-1}}$ is the standard metric in the unit $n-1$ sphere. Then
\[
\lam (\R^n, g_\a) \ge (n-1) \ln \a.
\]

(b).  Let $(M, g)$ be asymptotically $\e$ close to a round $\a$-cone. Then, there exists a constant $\beta>0$ depending only on $n$  such that
\[
\lam_\infty(M, g) \ge (n-1) \ln \a -\beta \e.
\]
\end{proposition}
\proof
Part (a).
 For any smooth compactly supported function $v$, with $\Vert v \Vert_{g_\a} =1$, the sharp log Sobolev functional can be written as
\be
\lab{logsobcone}
\al
&L(v, g_\a)= - \int_{\mathbf{R}^n} v^2 (\ln v^2) (\a r)^{n-1} dr dg_{S^{n-1}}\\
&+
\frac{n}{2} \ln \int_{\mathbf{R}^n} \left[ 4 \left( |\partial_r v|^2 + \frac{1}{(\a r)^2}
|\nabla_{g_{S^{n-1}}} v |^2 \right) + \frac{(n-2) (1-\a^2)}{(\a r)^2} v^2 \right] (\a r)^{n-1} dr dg_{S^{n-1}} +s_n.
\eal
\ee Here we have used the standard formula for the scalar curvature $R$ for warped product metrics. For example, using the formula for the Ricci curvature in Section 4.6 of \cite{CLN:1}, one sees that $R=(n-2) (1-\a^2)/(\a r)^2$. Consider the function
\[
v_0 = \a^{(n-1)/2} v.
\]Then
\[
\int_{\R^n} v^2_0 d g_{\R^n} = \int_{\R^n} v^2 \a^{n-1} r^{n-1} dr dg_{S^{n-1}} = \Vert v \Vert^2_{g_\a} =1.
\]Replacing $v$ in \eqref{logsobcone} by $v_0 \a^{-(n-1)/2}$, we find
\[
\al
&L(v, g_\a)= \ln \a^{n-1} - \int_{\R^n} v^2_0 (\ln v^2_0)  r^{n-1} dr dg_{S^{n-1}}\\
&+
\frac{n}{2} \ln \int_{\R^n} \left[ 4 \left( |\partial_r v_0|^2 + \frac{1}{(\a r)^2}
|\nabla_{g_{S^{n-1}}} v_0 |^2 \right) + \frac{(n-2) (1-\a^2)}{(\a r)^2} v^2_0 \right] r^{n-1} dr dg_{S^{n-1}} +s_n.
\eal
\]Since $\a \le 1$, we deduce
\[
L(v, g_\a) \ge L(v_0, g_{\R^n}) + (n-1) \ln \a \ge (n-1) \ln \a.
\]
This proves part (a) after taking the infimum over all smooth compactly supported $v$ with unit $L^2$ norm.
\medskip

Part (b).

 We first need to prove the following assertion.

{\it When the radius $\rho$ is sufficiently large, we have, for a dimensional constant $C$,
\be
\lab{assert2.1}
\lam(\M-B(0, \rho), g) \ge (1-C \e) \lam({\R}^n-J(B(0, \rho)), g_\a) - C \e.
\ee Here $J$ is the coordinate map near infinity in  Definition \ref{defAF}
of asymptotically $\e$ closeness to $\a$ round cones;  $g_\a$ is the round conic  metric.}

The proof is similar to \cite{Z14} Proposition 2.3 (b), (2.10).

Pick a function $v \in C^\infty_0(M-B(0, \rho))$ with $\Vert v
\Vert_{L^2(M-B(0, \rho), g)}=1$. Given $\e>0$, by Definition \ref{defAF} ,
for $x \in M-B(0, \rho)$ with $\rho$ sufficiently large, the following relations hold
\be
\lab{2.100}
(1-C \e) dg_\a \le dg(x) =\sqrt{ det g(x)} dx \le (1+C \e) dg_\a,
\ee Here and later in the proof, the constant $C$ depends only on $n$ unless stated otherwise.
\be
\lab{2.101}
(1-C \e) |\nabla_{g_\a} f | \le |\nabla_g v | \le (1+C\e)
|\nabla_{g_\a} f |
\ee
\be
\lab{2.101.1}
R_g(x) \ge  (1- C \e) R_{g_\a}(x) - C \e/d_{g_\a}(x, 0)^2, C \ge 1,
\ee
where $f = v \circ J^{-1}$ and $J$ is the coordinate map from $ M-B(0, r)$ to $\R^n$. Hence
\be
\lab{2.102}
\al
\int_{\M} &( 4 |\nabla_g v |^2 + R_g v^2 ) dg
\ge (1-C\e )^2
 \int_{\R^n} 4 |\nabla_{g_{\a}} f|^2 \sqrt{ det g(x)} dx\\
&+ (1-C \e)^2 \int_{\R^n} R_{g_\a} v^2 d g_\a  -C \e \int_{\R^n} \frac{v^2}{d_{g_\a}(x, 0)^2} \sqrt{ det g(x)} dx,
 \eal
\ee where $dx$ is the standard Euclidean volume element. Writing
\[
\sqrt{ det g(x)}/\sqrt{ det g_{\a}(x)} = w^2,
\] we deduce
\be
\lab{2.103}
\al
&\int_{\R^n} 4 |\nabla_{g_\a} f|^2 \sqrt{ det g(x)} dx = \int_{\R^n} 4 |w \nabla_{g_\a} f|^2  dg_\a\\
&=\int_{\R^n} 4 |\nabla_{g_\a}(w f) |^2 dg_\a - 8 \int_{\R^n} f \nabla_{g_\a} w \nabla_{g_\a} ( w f) dg_\a +
 4 \int_{\R^n} f^2 |\nabla_{g_\a} w |^2 dg_\a.
\eal
\ee By Definition \ref{defAF}, we know that
\[
|\nabla_{g_\a} w(x) | \le C \e/r
\]where $r=d_{g_\a}(x, 0)$. Hence,  we have
\be
\al
\lab{2.104}
\int_{\R^n} 4 |\nabla_{g_\a} f|^2 \sqrt{ det g(x)} dx &\ge (1-C \e) \int_{\R^n} 4 |\nabla_{g_\a}( f w) |^2 dg_\a - 16 \e^{-1} \int_{\R^n} f^2 |\nabla_{g_\a} w|^2 dg_\a,\\
&\ge (1-C \e) \int_{\R^n} 4 |\nabla_{g_\a}( f w) |^2 dg_\a - C \e^{-1} \e^2 \int_{\R^n} \frac{(f w)^2(x)}{d_{g_\a}(0, x)^{2}}  dg_\a.
\eal
\ee Note that the integrals are taking place outside of a large ball $B(0, \rho)$ under metric $g$ and $w$ are $\e$ close for $\rho$ large. For the same reason, the above integrals can be considered to take place outside the ball $B(0, \rho/2, g_\a)$ i.e. the one under the conic metric $g_\a$.

Using the same method as on the Euclidean case, i.e. integration by parts in the radial direction, one can prove the Hardy inequality on the round cone, which gives, since $f$ is compactly supported, that
\be
\lab{hardyincone}
\int_{\R^n} \frac{(f w)^2}{d_{g_\a}(0, x)^2}  dg_\a \le c_{n, \a} \int_{\R^n} |\nabla_{g_\a}( f w) |^2 dg_\a.
\ee Substituting this to the right hand side of \eqref{2.104}, we deduce
\[
\int_{\R^n} 4 |\nabla_{g_\a} f|^2 \sqrt{ det g(x)} dx \ge (1- C \e ) \int_{\R^n} 4 |\nabla_{g_\a}( f w) |^2 dg_\a.
\]This implies, after being plugged into \eqref{2.102}, that
\be
\lab{fv>1-e}
\int_{M} ( 4 |\nabla v |^2 + R v^2 ) dg  \ge (1-C\e)^2 (1-  C \epsilon ) \int_{\R^n} 4 |\nabla_{g_\a}( f w) |^2 dg_\a + (1-\e)^2 \int_{\R^n} R_{g_\a} d g_\a.
\ee where Hardy's inequality is again used on the last term, after writing $\sqrt{g(x)} dx = w^2 dg_\a$.
Also, since $w$ is $\e$ close to $1$ outside of large balls, $\ln w^2 \le C \e$,  we have
\be
\lab{2.105}
\al
\int_M v^2 \ln v^2 dg &=
 \int_{\R^n} (f w)^2 \ln f^2 dg_\a =\int_{\R^n} (f w)^2 \ln (f w)^2 dg_\a - \int_{\R^n} (f w)^2 \ln w^2 dg_\a \\
 &=\int_{\R^n} (f w)^2 \ln (f w)^2 dg_\a - O(1) \e.
 \eal
\ee This and (\ref{fv>1-e}) imply, after adjusting the constant $C$, that
 that
\be
\lab{2.106}
\al
&L(\M-B(0, \rho), v, g) =  -\int_M v^2 \ln v^2 dg + \frac{n}{2} \ln \int_{M} ( 4 |\nabla v |^2 + R v^2 ) dg  \\
&\ge -\underbrace{\int_{\R^n} (f w)^2 \ln (f w)^2 dg_\a}_{N(fw)} + \frac{n}{2} \ln  \underbrace{\left[\int_{\R^n} \left(4 |\nabla_{g_\a}( f w) |^2 + R_{g_\a} \right)dg_\a \right]}_{F(fw)} -C \e\\
&=-N(fw) + F(fw)  - C \e\\
&= L({\R}^n-J(B(0, \rho)), f w, g_\a)  - C \e.
\eal
\ee
Since $\Vert f w \Vert_{L^2(\R^n, g_\a)} = 1$, by taking the infimum of  inequality \eqref{2.106},
we find
\be
\lab{2.107}
\lam(\M-B(0, \rho), g) \ge  \lam(
{\R}^n-J(B(0, \rho)), g_\a,)  - C \e.
\ee  The assertion is proven.

Using
\[
\lam(g_\a,  {\R}^n-J(B(0, \rho)) \ge \lam(g_a,  {\R}^n) \ge (n-1) \ln \a,
\]which is due to part (a), we see, after renaming the constant, that
\be
\lab{2.108}
\lam_\infty(g) = \lim_{\rho \to \infty} \lam(g,  \M-B(0, \rho)) \ge (n-1) \ln \a - \beta \e.
\ee
This proves the proposition.
 \qed

Next let us recall a previous result on existence of minimizers  which is needed  later in this paper.

\begin{theorem}  (Theorem 1.10 \cite{Z14})
\lab{thminimizer} Let $(M, g)$ be a noncompact manifold with
bounded curvature and nonnegative scalar curvature, which also
satisfies

 (a) $-\infty<\lam(M, g) < \lam_\infty(M, g)$.

 (b) Either $|B(x_0, r)|_{g} \le C r^n$, for some $C>0$ and all $r>0$,
 or  $R(x) \ge \frac{C}{1+ d(x, x_0)^2}$ for some constant $C>0$.

Then there exists a smooth minimizer $v$ for the Log Sobolev functional $L(M, \cdot, g)$ in $W^{1, 2}(M)$, which satisfies
the equation
\be
\lab{maineq}
 \frac{n}{2} \frac{4 \Delta v - R v}{ \int ( 4 | \nabla v |^2 + R v^2 ) dg }
+  2 v \ln v +   \left( \lam(M, g)+   \frac{n}{2} -  \frac{n}{2}  \ln  \int ( 4 | \nabla v |^2 + R v^2 ) dg
- s_n \right)  v = 0.
\ee
\end{theorem}

In addition, we also need an extension of  the above theorem to a suitable family of manifolds so that the maximum of values of some minimizers have uniform positive lower bound.

\begin{theorem}
\lab{thminimizerk} Let $(M, g)$ be a noncompact n-manifold with the following properties. For some positive constants $c_1, ..., c_5$,

(a)  $|\nabla^k Rm| \le c_1$ for $k=0, 1, 2, 3$;

(b) $|B(x, 1)| \ge c_2$ for all $x \in M$;

(c) $|B(x, r)| \le c_3 r^n$ for all $x \in M$, $r>0$;

(d) There is a point $x_0 \in M$ such that
\[
-c_4<\lam(M, g) < \lam(M-B(x_0, r_*), g)-c_5<0
\]for some fixed number $r_*>1$

Then there exists a minimizer $v$ for the Sharp log Sobolev functional $L(M, v, g)$ such that
\[
\sup v \ge m>0
\] where the constant $m$ depends only on $r_*$ and $c_i$, $i=1, ..., 5$.
\proof
\end{theorem}

Suppose the result is not true. Then we can find a sequence of manifolds $(M_k, g_k)$ satisfying all the assumptions in the theorem, together with minimizers of the sharp log functional $v_k$ such that
\be
\lab{vkto0}
\lim_{k \to \infty}  \Vert v_k \Vert_\infty =0.
\ee

{\it Step 1.}

 Since $v_k$ is a minimizer for $L(M_k, v,  g_k)$, we have
 \be
 \lab{lamk=1}
 \al
 \lam_k& \equiv \lam_k(M_k, g_k)=L(M_k, v_k,  g_k)
\\
&= - \int_{M_k} v^2_k \ln v^2_k dg_k + \frac{n}{2} \ln
\left( \int_{M_k} (4 |\nabla v_k|^2 + R_k v^2_k ) dg_k \right) + s_n.
 \eal
 \ee
According to Theorem \ref{thminimizer}, $v_k (>0)$ exists. Using Lemma 3.3 in \cite{Z14},  we know that $\{ \Vert v_k \Vert_\infty \}$ is uniformly
bounded. We comment that the cited lemma was stated for balls in a fixed manifold but the proof is almost identical under the conditions of the current theorem. The only change is to take $C^2_{loc}$ limits for $M_k$ instead of for the balls.

From \eqref{vkto0},   there exists a sequence of positive
integers $\{ i_k \}$ and a subsequence of $\{ v_k \}$,  denoted by the same
symbol,
such that $i_k \to \infty$ slow enough as $k \to \infty$
and that
\be
\lab{vkBto0} \int_{B(x_k,  2^{2 i_k}, g_k)}  v^2_k dg_k \to
0, \qquad k \to \infty.
\ee Here $x_k$ is the  reference point in place of $x_0$ in condition (d) for the manifold $M_k$.
 Indeed, as long as $|B(x_k,  2^{2 i_k}, g_k)|_{g_k}=o(1/\Vert v_k \Vert_\infty)$, then the above limit holds.
 For any positive integer $i$ we
introduce the following notations
\be
\lab{3.141}
\al
\Omega_{ki} &=B(x_0, 2^i, g_k) - B(x_0, 2^{i-1}, g_k),\\
F(v_k) & = \int_{M_k} ( 4 |\nabla v_k |^2 + R_k v_k^2 ) dg_k, \quad
N(v_k) =  \int_{M_k} v^2_k \ln v^2_k dg_k.
\eal
\ee

By $\lam_k \ge -c_4$ in assumption (d) of
the theorem and Proposition 2.6 in \cite{Z14}, there exists a positive constant $A=A(c_1, c_4, n)$ such that
\be
\lab{3.142}
\left( \int_{M_k} v^{2n/(n-2)}_k dg_k \right)^{(n-2)/n} \le
A F(v_k).
\ee We comment that, this statement is just a variation of, involving the scalar curvature, of the well known fact that sharp log Sobolev inequality implies the standard $L^2$ Sobolev inequality.

Hence \be \lab{Fv<eNv} \al &\left( \Sigma^{2 i_k}_{i=i_k}
\int_{\Omega_{ki}} v^{2n/(n-2)}_k dg_k \right)^{(n-2)/n} e^{-N(v_k)
2/n}   \\
&\le A F(v_k) e^{-N(v_k) 2/n} = A e^{(\lam_k -s_n) 2/n} \le C=C(c_1, c_4, n),
 \eal
\ee where, in the last line, we also used (\ref{lamk=1}) and the fact that $\lam_k \le 0$. Thus, there exists an integer $j_k \in [i_k, 2
i_k]$ such that
\be \lab{vkojk<} \left(  \int_{\Omega_{kj_k}}
v^{2n/(n-2)}_k dg \right)^{(n-2)/n} \le C i_k^{-(n-2)/n}
e^{N(v_k) 2/n} \ee

By partition of unity, we can choose a sequence of cut-off
functions $\phi_k$, $\eta_k$ on $M_k$ such that $ \phi_k = 1 $ on
$B(x_0, 2^{j_k -1}, g_k)$,  \, $supp \, \phi_k \subset B(x_0,
2^{j_k}, g_k)$; $ \eta_k = 1 $ on $M-B(x_0, 2^{j_k}, g_k)$,  \, $supp \,
\eta_k \subset M_k- B(x_0,  2^{j_k-1}, g_k)$;
 $| \nabla_{g_k} \phi_k | + | \nabla_{g_k} \phi_k | \le C/{2^{j_k}}$;
$\phi^2_k + \eta^2_k=1$. Here $C$ is an absolute constant.  We introduce the notations
\be
\lab{3.143}
a_k \equiv \Vert v_k \phi_k \Vert^2_{L^2(g_k)}, \quad b_k \equiv \Vert
v_k \eta_k \Vert^2_{L^2(g_k)};
\ee
\be
\lab{3.144}
A_k \equiv \exp( \frac{2}{n} N(v_k \phi_k)), \quad B_k \equiv
\exp( \frac{2}{n} N(v_k \eta_k)).
\ee By (\ref{vkBto0}) and $\Vert v_k \Vert_{L^2(g_k)} =1$, we know that
\be
 \lab{akbkto01} a_k \to 0, \quad b_k \to 1, \quad \text{as}
\quad k \to \infty.
\ee

Now we will split the terms in the sharp log
Sobolev functional into terms involving $v_k \phi_k$ and $v_k
\eta_k$. By direct computation
\be
\lab{fenF}
\al
&\int (4 | \nabla_{g_k} v_k |^2 + R_k v^2_k ) dg_k \\
&=\int (4 | \nabla_{g_k} (v_k \phi_k) |^2 + R_k (v_k \phi_k)^2 ) dg_k +
\int (4 | \nabla_{g_k} (v_k \eta_k) |^2 + R_k (v_k \eta_k)^2 ) dg_k\\
&\qquad - 4 \int (  | \nabla_{g_k} \phi_k|^2 + | \nabla_{g_k} \eta_k|^2 )  v^2_k dg_k,
\eal
\ee where we have used the identity
\be
\lab{3.145}
0=\Delta_{g_k} (\phi^2_k+\eta^2_k) = 2  | \nabla_{g_k} \phi_k|^2  + 2 \phi_k \Delta_{g_k} \phi_k + 2 | \nabla_{g_k} \eta_k|^2
+ 2 \eta_k \Delta_{g_k} \eta_k.
\ee

Using Condition (c):
$|B(x_k, r, g_k)| \le c_3 r^n$ and H\"older's inequality we deduce
\be
\lab{3.146}
\al
&4 \int (  | \nabla_{g_k} \phi_k|^2 + | \nabla_{g_k} \eta_k|^2 )  v^2_k dg_k \le
C 2^{-2 j_k} \int_{\Omega_{kj_k}} v^2_k dg_k\\
&\le C 2^{-2 j_k} |\Omega_{kj_k}|^{2/n} \left(  \int_{\Omega_{kj_k}} v^{2n/(n-2)}_k dg_k \right)^{(n-2)/n} \le C  \left(  \int_{\Omega_{kj_k}} v^{2n/(n-2)}_k dg_k \right)^{(n-2)/n}.
\eal.
\ee  Using (\ref{vkojk<}), we know that
\be
\lab{dphi+deta} 4 \int (  | \nabla \phi_k|^2 + | \nabla
\eta_k|^2 ) v^2_k dg \le  C i_k^{-(n-2)/n}
 e^{N(v_k) 2/n}.
\ee Here $o(1)$ is a quantity that goes to $0$ when $k \to
\infty$. This and (\ref{fenF}) imply
\be \lab{splitFvk} F(v_k) =
F(v_k \phi_k)  + F(v_k \eta_k) - \delta_k
 e^{N(v_k) 2/n}, \qquad  0 \le \delta_k \le  C i_k^{-(n-2)/n}.
\ee

Next,
observe, since $\Vert v_k \Vert_\infty \le 1$ by \eqref{vkto0} for $k$ large, that
\be
\lab{3.149}
\al
\int v^2_k &\ln v^2_k dg_k -
\int (v_k \phi_k)^2 \ln (v_k \phi_k)^2 dg_k - \int (v_k \eta_k)^2 \ln (v_k \eta_k)^2dg_k\\
&=\int (v_k \phi_k)^2 \left[ \ln ( (v_k \phi_k)^2 +
(v_k \eta_k)^2 ) -  \ln (v_k \phi_k)^2 \right] dg_k  \\
&\qquad + \int (v_k \eta_k)^2 \left[ \ln ( (v_k \phi_k)^2  +
(v_k \eta_k)^2 ) -  \ln (v_k \eta_k)^2 \right] dg_k \\
&\le 2 \int v^4_k \phi^2_k \eta^2_k dg_k \le 2 \Vert v_k \Vert_\infty^2 \int_{\Omega_{j_k}} v^2_k \, dg_k \le 2 \Vert v_k \Vert_\infty^2. \eal
\ee Here we just used the inequality $\ln(p+q)-\ln p \le q$ for $p \in (0, 1]$ and $q>0$.
This means
\be
\lab{splitNvk}
N(v_k) = N(v_k \phi_k) + N(v_k \eta_k) + \e_k, \qquad 0 \le \e_k \le 2 \Vert v_k \Vert_\infty^2.
\ee

Recall that $v_k$ is a minimizer for the log Sobolev functional. By (\ref{lamk=1}),
\be
\lab{Fv/eN}
e^{\frac{2}{n} (\lam_k -s_n)} = \frac{F(v_k)}{\exp (\frac{2}{n}
N(v_k)) }.
\ee By (\ref{splitFvk}) and (\ref{splitNvk}), this implies
\be
\lab{3.150}
\al
e^{\frac{2}{n} (\lam_k-s_n)}&=\frac{F(v_k \phi_k) + F(v_k \eta_k)
-\delta_k \exp( \frac{2}{n} N(v_k))}{\exp( \frac{2}{n} N(v_k))} \\
&=\frac{F(v_k \phi_k) + F(v_k \eta_k)
 }{\exp( \frac{2}{n} N(v_k \phi_k)) \, \exp( \frac{2}{n}
N(v_k \eta_k)) \, e^{\e_k}} -\delta_k.
\eal
\ee  On the other hand, by definition of $\lam_k$, we have
\be
\lab{3.151}
F(v_k \phi_k) \ge e^{\frac{2}{n} (\lam_k-s_n)} \Vert v_k \phi_k
\Vert^2_{L^2(g_k)} \exp \left( -\frac{2}{n} \ln \Vert v_k \phi_k
\Vert^2_{L^2(g_k)} \right) \exp\left(\frac{2}{n} N(v_k \phi_k) / \Vert
v_k \phi_k \Vert^2_{L^2(g_k)}\right).
\ee
Write
\[
\lam_{k,r_*, \infty}=\lam(M_k-B(x_k, r_*, g_k), g_k).
\]
Since the support of $\eta_k$ is outside of the ball $B(x_0, 2^{j_k-1}, g_k)$, by Definition \ref{deflv}, we
know
\be
\lab{3.152}
F(v_k \eta_k) \ge e^{\frac{2}{n} (\lam_{k,r_*, \infty}-s_n)} \Vert v_k \eta_k
\Vert^2_{L^2(g_k)} \exp \left( -\frac{2}{n} \ln \Vert v_k \eta_k
\Vert^2_{L^2(g_k)} \right) \exp\left(\frac{2}{n} N(v_k \eta_k) / \Vert
v_k \eta_k \Vert^2_{L^2(g_k)}\right).
\ee  Combining \eqref{3.150}, \eqref{3.151}, and \eqref{3.152}, we deduce, that
\be
\lab{3.153}
1 \ge \frac{a^{-2/n}_k a_k A^{1/a_k}_k + b^{-2/n}_k b_k B^{1/b_k}_k
e^{(\lam_{k, r_*, \infty}-\lam_k) 2/n} }{A_k B_k e^{\e_k}} -\delta_k,
\ee  where
\be
\lab{3.154}
a_k \equiv \Vert v_k \phi_k \Vert^2_{L^2(g_k)}, \quad b_k \equiv \Vert v_k \eta_k \Vert^2_{L^2(g_k)};
\ee
\be
\lab{3.155}
A_k \equiv \exp( \frac{2}{n} N(v_k \phi_k)), \quad B_k \equiv \exp( \frac{2}{n} N(v_k \eta_k)).
\ee  Therefore, by part of condition (d), i.e.
\[
\lam_{k, r_*, \infty}-\lam_k \ge c_5>0,
\]we deduce
\be
\lab{3.156}
\min\{ a^{-2/n}_k, \,b^{-2/n}_k \} \,
\frac{ a_k A^{1/a_k}_k +  b_k B^{1/b_k}_k e^{ c_5 2/n} }{A_k B_k e^{\e_k}}  - \delta_k \le 1,
\ee Since $a_k$ and $b_k$ are positive numbers in the interval $(0, 1)$, this shows
\be
\lab{3.157}
\ln (a_k A^{1/a_k}_k +  b_k B^{1/b_k}_k e^{c_5 2/n}) \le \ln (A_k B_k) + \delta_k+\e_k.
\ee
 Notice that $a_k + b_k =1$. By concavity of the $\ln$ function we obtain
\be
\lab{3.158}
b_k c_5 2/n  \le \delta_k+\e_k.
\ee
Letting $k \to \infty$ and using the fact that $b_k \to 1$ (from (\ref{akbkto01}) ) and $\e_k, \delta_k \to 0$,  we arrive at
\be
\lab{3.159}
c_5 \le 0.
\ee Since $c_5$ is a positive number by assumption,  this is a contradiction  which proves that \eqref{vkto0} is false. The conclusion of the theorem is true. \qed

\medskip

Now we are ready to give

\noindent {\bf  Proof of Theorem \ref{thmrig}}.

Suppose the conclusion of of the theorem is false. Then there is a non-flat ancient solution $(M, g(t))$ satisfying the conditions of the theorem.  Let us choose a sequence of time $t_k$ for the ancient solution, which goes to $-\infty$ as $k \to \infty$.

Choosing $\e_0=\delta_*/2$ in our assumption, we deduce, from Proposition \ref{prAFsob} (b), that
\be
\lab{lam-lam0}
\lam(M-B(0, r_0), g(t_k)) - \lam(M, g(t_k)) \ge \delta_*/2.
\ee

 This inequality allows us to use Theorem \ref{thminimizerk}, which says that a minimizer $v_k$ for $\lam_k \equiv \lam(M, g(t_k))$ satisfies, for a fixed number $m_1$,
\be
\lab{vk>m1}
\sup_k  \Vert v_k \Vert_\infty \ge m_1>0
\ee Since $(M, g(t_k))$ has bounded geometry, it is not hard to see, c.f. Lemma 3.3. \cite{Z14}, that there is another constant $m_2>0$ such that
\be
\lab{vk<m2}
\sup_k  \Vert v_k \Vert_\infty \le m_2.
\ee Note that, from the assumption that the flow is $\kappa$ non-collapsed at scale 1 and the condition \eqref{lambbce0} near infinity, we can deduce that
\be
\lab{lammg>-C}
\lam(M, g(t)) \ge -C_2>-\infty, \quad \forall t \le 0.
\ee This follows from a standard method of splitting $M$ into the union of compact and noncompact domains as in the middle of the proof of Theorem
\ref{thminimizerk}.

Consider the shifted metrics
\[
\widetilde{g}_k=\widetilde{g}_k(\cdot, t)=g(\cdot, t_k + t), \quad t \le 0,
\]which has bounded geometry too. Using Hamilton's compactness theorem, there exists a subsequence, still denoted by $\{ \widetilde{g}_k \}$ such that the pointed sequence of Ricci flows $\{(M, \widetilde{g}_k, 0)\}$, converges in $C^\infty_{loc}$ topology, to a limit ancient Ricci flow
\be
\lab{mwuq}
(M_\infty, g_\infty, p_\infty),
\ee which  will next be shown to be a gradient Ricci soliton.
The proof of this assertion is similar to the proof of part (a) of the theorem.

We solve the conjugate heat equation
\[
\Delta_{\widetilde{g}_k(t)} u - R_{\widetilde{g}_k(t)} u + \partial_t u =0
\]
 for $t \le 0$, with final value as $u_k(\cdot, 0)=v^2_k(\cdot)$. This solution
is denoted by $u_k=u_k(x, t)$.
Write $v_k = \sqrt{u_k}$, then by Definition \ref{deflv}
\be
\lab{3.179.1}
L(v_k, \widetilde{g}_k(t)) = - N(v_k) + \frac{n}{2} \ln F(v_k) + s_n,
\ee  where, due to $v_k = \sqrt{u_k}$,
\be
\lab{3.180.1}
\al
&N(v_k) = \int_M u_k \ln u_k \, d\widetilde{g}_k(t); \quad\\
 &F(v_k) = \int_M (\frac{|\nabla_{\widetilde{g}_k} u_k|^2}{u_k} + R_{\widetilde{g}_k} u_k) d\widetilde{g}_k(t) =
\int_M (4 |\nabla_{\widetilde{g_k}} v_k|^2 + R_{\widetilde{g_k}} v^2_k) d\widetilde{g}_k(t).
\eal
\ee  Then
\be
\lab{3.183.1}
\al
\int^{0}_{t_{k+1}-t_k} \frac{d}{dt}  L(\sqrt{u_k}, \widetilde{g}_k(t)) dt &=
L(\sqrt{u_k(\cdot, 0)}, \widetilde{g}_k(0)) - L(\sqrt{u_k(\cdot, t_k-t_{k+1})}, \widetilde{g}_k(t_k-t_{k+1}))\\
& \le \lam(M, \widetilde{g}_k(0)) - \lam(M, \widetilde{g}_k(t_k-t_{k+1}))\\
& = \lam(M, g(t_k)) - \lam(M, g(t_{k+1})) \to 0,  \quad k \to \infty.
\eal
\ee
According to Perelman \cite{P:1} Section 1, $\frac{d}{dt} N(v_k) = F(v_k)$ and
\be
\lab{3.181.1}
\frac{d}{dt} F(v_k) = 2 \int_M | Ric_{\widetilde{g_k}} - Hess_{\widetilde{g_k}} (\ln u_k)|^2 u d\widetilde{g_k}(t).
\ee
 Following Perelman's calculation in \cite{P:1}, we arrive at
\be
\lab{3.181.2}
\al
0 \le
&\int^0_{t_{k+1}-t_k} \left[ n \int_M | Ric_{\widetilde{g_k}} - Hess_{\widetilde{g_k}}(\ln u_k) - \frac{1}{n} ( R_{\widetilde{g_k}} - \Delta_{\widetilde{g_k}} \ln u_k) g_{\widetilde{g_k}} |^2
u_k dg_{\widetilde{g_k}}(t) \right] F^{-1}(v_k) dt\\
&
+\int^0_{t_{k+1}-t_k} \left[ \int_M ( R_{\widetilde{g_k}} - \Delta_{\widetilde{g_k}} \ln u_k)^2 u_k \, d{\widetilde{g_k}}(t)  -  \left(\int_M ( R_{\widetilde{g_k}} - \Delta_{\widetilde{g_k}} \ln u_k) u_k \, d{\widetilde{g_k}}(t) \right)^2\right]  F^{-1}(v_k) dt\\
&\le \int^{0}_{t_{k+1}-t_k} \frac{d}{dt}  L(\sqrt{u_k}, \widetilde{g}_k(t)) dt \to 0,
\eal
\ee where we just used \eqref{3.183.1}.  It is well known that although Perelman only proved the formulas for compact manifolds,  his proof also works for noncompact manifolds with bounded geometry when the functions involved have sufficiently fast decay such as the Gaussian function. See \cite{chowetc3} Chapter 19 and
 \cite{cty2007}
e.g.. In our case, the functions involved $u_k$ have Gaussian type decay at each time level.

Due to \eqref{vk>m1} and \eqref{vk<m2} and the equation \eqref{maineq} for $v_k(\cdot, 0)$ as a minimizer,
by standard elliptic theory,  a subsequence of $\{ v_k
\}$, still denoted by the same symbol, converges in
$C^\infty_{loc}$ sense, to a nontrivial limit function $v_\infty \in
C^\infty(M_\infty)$. Since $u_k=v^2_k$, writing $u_\infty=v^2_\infty$, from \eqref{3.181.2}, we deduce
\be
\lab{3.181.3}
\al
&\int^0_{-\infty} \left[ n \int_{M_\infty} | Ric_{g_\infty} - Hess_{g_\infty}(\ln u_\infty) - \frac{1}{n} ( R_{g_\infty} - \Delta_{g_\infty} \ln u_\infty) g_\infty |^2
u_\infty dg_\infty(t) \right] F^{-1}(v_\infty) dt\\
&
+\int^0_{-\infty} \left[ \int_{M_\infty} ( R_{g_\infty} - \Delta_{g_\infty} \ln u_\infty)^2 u_\infty \, dg_\infty(t)  -  \left(\int_{M_\infty} ( R_{g_\infty} - \Delta_{g_\infty} \ln u_\infty) u_\infty \, dg_\infty(t) \right)^2\right]  \\
&\qquad \qquad F^{-1}(v_\infty) dt\\
&= 0.
\eal
\ee Here we have selected a subsequence of $\{t_k\}$, without changing the notation, so that $t_k-t_{k+1} \to \infty$. Some justification of the above limit process is needed.
The first integral of \eqref{3.181.2} is taken care of by Fatou's lemma. For the second one, since the overall integrand in the time integral is nonnegative, we can take the limit sign inside the time integral by Fatou's lemma with a correct inequality. For the space integrals, there is a negative sign inside. In order to take the limit, we notice that for each fixed time, $v_k$ decays exponentially near infinity. This is due to the exponential  decay of the final value $v_k(\cdot, 0)$ as a minimizer (c.f. \cite{Z:2} e.g.) and finite time exponential decaying of the conjugate heat kernel on manifold with bounded geometry.  Notice that we do not need uniform spatial decay for all time. Due to equality case of the Cauchy Schwarz inequality, we find that
\be
\lab{rddlnu0}
R_\infty- \Delta_{g_\infty} (\ln u_\infty) -  l(t) =0.
\ee where $l(t)$ is a function of time only.

Therefore, we arrive at
\[
Ric_{g_\infty} - Hess_{g_\infty}(\ln u_\infty) - \frac{1}{n}l(t) g_\infty=0, \quad t \le 0.
\]
 Since $u_\infty$ has Gaussian type decay near infinity, we can multiply both sides of \eqref{rddlnu0} by $u_\infty$ and integrate, giving us
\[
\int_{M_\infty}  \left(R_\infty u_\infty + \frac{|\nabla_{g_\infty} u_\infty|^2}{u_\infty} \right) dg_\infty(t)=  l(t).
\]Thus $l(t)>0$ and $(M_\infty, g_\infty(t))$ is a gradient shrinking soliton.

 In particular, since $t$ extends to $-\infty$, we know that
\[
\lim_{t \to -\infty} |Rm_{g_\infty}(p_\infty, t)|=0.
\]Let us recall that the time $t$ for $g_\infty$ is actually the limit of the shifted time $t+t_k$ for the original metric $g$ of the ancient solution, with $t_k \to -\infty$. Also $p_\infty$ corresponds to the reference point $0$ in $M$. Therefore, there exists a sequence of time $s_k$ going to $-\infty$ such that
\be
\lab{rm0to0}
\lim_{k \to \infty} |Rm_{g}(0, s_k)|=0.
\ee Using these $s_k$ as the final times and consider the shifted ancient solutions
\[
\widetilde{g_k}=g(\cdot, t+s_k), \quad t \le 0.
\]Repeating the above process, we can find another subsequence, still denoted by $(M, \widetilde{g_k}, 0)$, which converges to a gradient shrinking soliton in $C^\infty_{loc}$ sense. We still denote this limit as $(g_\infty, p_\infty)$. By \eqref{rm0to0}, $|Rm_{g_\infty}(p_\infty, 0)|=0.$ In particular the scalar curvature is $0$ at that point. Since the scalar curvature is nonnegative, by the maximum principle, applied on the scalar curvature equation, we deduce that the scalar curvature for $g_\infty$ is $0$ everywhere. Using the equation for the scalar curvature along a limiting Ricci flow:
\[
\Delta_{g_\infty} R_{g_\infty} - \partial_t R_{g_\infty} + 2 |Ric_{g_\infty}|^2=0,
\]we infer that the Ricci curvature for $g_\infty$ is $0$.
Fixing a time $t$, we find, from defining equation for gradient shrinking soliton
\[
 Hess f =
c g_\infty
\]for some smooth $f$ and positive constant $c$. Applying Theorem 2 (I b) in \cite{T:1} , we conclude that $g_\infty$ is the Euclidean metric and $M$ is $\R^n$.

Let $v_k$ be a minimizer for the sharp log functional on $(M, g(s_k))$. Then from
\eqref{lamk=1}, we know that
\be
\lab{lamk=11}
 \al
 \lam_k& \equiv \lam_k(M, g(s_k))=L(M, v_k, g(s_k))
\\
&= - \int_{M} v^2_k \ln v^2_k dg(s_k) + \frac{n}{2} \ln
\left( \int_{M} (4 |\nabla v_k|^2 + R v^2_k ) dg(s_k) \right) + s_n.
\eal
\ee As before, $v_k$ sub-converges in $C^\infty_{loc}$ sense to a nontrivial function $v_\infty$ on $(M_\infty, g_\infty(0))$. Due to its exponential decaying near infinity, we can take limit in \eqref{lamk=11} to reach
\be
\lab{lamk=12}
 \al
\lim_{k \to \infty} \lam_k = - \int_{\R^n} v^2_\infty \ln v^2_\infty dx + \frac{n}{2} \ln
\left( \int_{\R^n} (4 |\nabla v_\infty|^2 ) dx \right) + s_n \ge 0.
\eal
\ee Here we just used the sharp log Sobolev inequality on $\R^n$.
However, by Proposition
\ref{prAFsob} (b), we know $\lim_{k \to \infty} \lam_k \le -\delta_* <0$. This is a contradiction which shows that $(M, g(t))$ is flat.  From \eqref{lammg>-C} and flatness, the standard Euclidean Sobolev inequality holds on $(M, g(t))$ which then has maximum volume growth. Therefore $M$ is the standard $\mathbf{R}^n$ due to flatness again.

The converse is obviously true since $\lam(\mathbf{R}^n)=0 \ge -\e_0$. The proof is complete.
\qed

\medskip

Next we give

\noindent {\bf  Proof of Corollary \ref{corig} (a).}

This is similar to the proof of the theorem.

We use the method of contradiction. Suppose, under the condition of the theorem, $(M, g)$ are  gradient Ricci solitons. Then we can solve the Ricci flow with $g$ as the initial metric to form a smooth Ricci flow $(M, g(t))$ in an open time interval.

{\it Step  1.   existence of a minimizer}

First we recall  that $\lam(g)$ is invariant under scaling and  diffeomorphism. The proof is quite easy and can be found in \cite{Z14}  Section 3 near (3.118), e.g.
 Hence, we know from the contrapositive assumption of $(M, g(t))$ being solitons:
$g(t_2) = c \psi^* g(t_1)$  that
\be
\lab{tk-tk}
 \lam(g(t_1)) - \lam(g(t_2))=0.
\ee Here $c$ is a positive constant, $t_2>t_1$ and $\psi$ is a diffeomorphism.

Since $(M, g)$ is a fixed smooth manifold, we know that it is $\kappa$ non-collapsed at scale $1$ in the ball $B(0, 2r_0)$ for some $\kappa>0$. Here $r_0$ is from Definition \ref{defAF}. Moreover, $(B^c(0, r_0), g)$ is $\e$ close to a $\alpha$ cone by assumption. These imply that $\lam(B(0, 2r_0), g) \ge -C_1$ and $\lam(B^c(0, r_0), g) \ge -C_1$ for some positive constant $C_1$. Then a standard splitting argument show that
\[
\lam(M, g) \ge -C_2>-\infty.
\]See the argument in the middle of the proof of Theorem \ref{thminimizerk} e.g. Our assumption on the $\e$ closeness to a $\alpha$ cone near infinity, Proposition \ref{prlamcaaca} and Proposition \ref{prAFsob} (b), (c) imply that
\[
\lam(M, g) < \lam_\infty(M, g).
\]
Then according to Theorem \ref{thminimizer}, there exists a function $v_2 \in W^{1, 2}(M, g(t_2))$,
which is a minimizer for $\lam(g(t_2))$, i.e.
\be
\lab{l=lam}
L(v_2, g(t_2))  = \lam(g(t_2)).
\ee Moreover, by Moser's iteration, it is known, as done in Lemma 2.3 in \cite{Z:2}, that $v_2$ has Gaussian type decay at infinity.

Next, we solve the conjugate heat equation
\[
\Delta u - R u + \partial_t u =0
\]
 for $t<t_2$, with final value as $v^2_2$. This solution
is denoted by $u=u(x, t)$.
Write $v = \sqrt{u}$, then by Definition \ref{deflv}
\be
\lab{3.179}
L(v, g(t)) = - N(v) + \frac{n}{2} \ln F(v) + s_n,
\ee  where, due to $v = \sqrt{u}$,
\be
\lab{3.180}
N(v) = \int_M u \ln u \, dg(t); \quad F(v) = \int_M (\frac{|\nabla u|^2}u + R u) dg(t) =
\int_M (4 |\nabla v|^2 + R v^2) dg(t).
\ee  According to Perelman \cite{P:1} Section 1, $\frac{d}{dt} N(v) = F(v)$ and
\be
\lab{3.181}
\frac{d}{dt} F(v) = 2 \int_M | Ric - Hess (\ln u)|^2 u dg(t).
\ee As in the proof of the theorem, in our case, the function $v$ has Gaussian type decay at each time level just like the final value $v(t_2)$ does. Hence all integrations are rigorous and
\be
\lab{d/dtlv=}
\frac{d}{dt} L(v, g(t)) = \left( n \int_M | Ric - Hess (\ln u)|^2 u dg(t) - F^2(v) \right) \, F^{-1}(v).
\ee Following Perelman's computation again,
\be
\lab{3.182}
| Ric - Hess (\ln u)|^2  \ge \left| Ric - Hess (\ln u) - \frac{1}n ( R - \Delta \ln u) g \right|^2
+ \frac{1}n ( R - \Delta \ln u)^2;
\ee Using the relation $F(v) = \int_M  (R - \Delta \ln u) u \, dg(t)$, we deduce
\be
\lab{dldt>}
\frac{d}{dt} L(\sqrt{u}, g(t)) \ge \frac{Q(u)}{F(v)}  \ge 0
\ee where
\be
\lab{qu}
\al
Q(u)(t) &= n \int_M | Ric - Hess (\ln u) - \frac{1}{n} ( R - \Delta \ln u) g |^2
u dg(t) \\
&\qquad
+ \int_M ( R - \Delta \ln u)^2 u \, dg(t)  -  \left( \int_M  (R - \Delta \ln u) u \, dg(t) \right)^2;\\
F(v)&=F(v)(t)=F(\sqrt{u})(t) = \int_M (\frac{|\nabla u|^2}u + R u) dg(t).
\eal
\ee Observe that $\sqrt{u(\cdot, t_2)} = v_2(\cdot)$ by definition. So by (\ref{l=lam}) we deduce
\be
\lab{3.183}
\al
\int^{t_2}_{t_1} \frac{d}{dt}  L(\sqrt{u}, g(t)) dt &=
L(\sqrt{u(\cdot, t_2)}, g(t_2)) - L(\sqrt{u(\cdot, t_1)}, g(t_1))\\
& \le \lam(g(t_2)) - \lam(g(t_1))=0.
\eal
\ee The last line is due to  (\ref{tk-tk}). By (\ref{dldt>}), we then have
\be
\lab{q/f}
F^{-1}(v) Q(u) =0.
\ee By (\ref{qu}) and the equality case of the Cauchy Schwarz inequality, this shows that
\[
(R - \Delta \ln u)(\cdot, t)=l(t)
\] where $l=l(t)$ is a function of $t$ only.
In addition, it implies
\be
\lab{3.184}
Ric - Hess (\ln u) - \frac{1}{n} l(t) g=0.
\ee Therefore, $(M, g(t))$ is a gradient Ricci soliton.  Taking the trace, we deduce
\be
\lab{rddlnu}
R- \Delta (\ln u) - l(t) =0.
\ee Since $u$ has Gaussian type decay near infinity, we can multiply both sides of \eqref{rddlnu} by $u$ and integrate, giving us
\[
\int_M  (R u + \frac{|\nabla u|^2}{u}) dg(t)= l(t).
\]Thus $l(t)>0$ and $(M, g(t))$ is a gradient shrinking soliton.
\medskip

{\it Step  2.  ruling out non-flat gradient steady and expanding solitons.}

 From now on, we fix a time level $t_0$ and suppress the time variable. Since we also assumed $M$ is a steady or expanding soliton, there exists a potential function $f$ such that
\be
\lab{richessfc}
Ric+ Hess f =0, \quad \text{or} \quad Ric+ Hess f +
c g =0, \quad c>0.
\ee Subtracting this by \eqref{3.184}, we see that the Ricci curvature is cancelled and that
\be
\lab{hessfcg}
Hess \, F = c_1 g, \quad F=-(f+ \ln u),
\ee where $c_1$ is another positive constant.
 According to Theorem 2 (I b) in \cite{T:1} by Tashiro again, the manifold $M$ is isometric to $\R^n$ with the standard metric. This is a contradiction with the non-flatness assumption. Hence $M$ can not be a non-flat gradient steady or expanding Ricci soliton.
\qed


\medskip

\noindent {\bf  Proof of Corollary \ref{corig} (b).}
\medskip
This is a direct consequence of the theorem and the propositions at the beginning of the section.

According to Proposition \ref{prAFsob} (b),
\be
\lim_{k \to \infty} \lam(M, g(t_k))
\le -\delta_*.
\ee By our assumption and Proposition \ref{prlamcaaca} (b),
\be
\lam(M-B(0, r_*), g(t_k)) \ge  (n-1) \ln \a -\beta \e.
\ee Therefore, when $\e$ is sufficiently small and $\a>\a_0$ which is sufficiently close to $1$, we have
\be
\lab{lam-lam}
\lam(M-B(0, r_*), g(t_k)) - \lam(M, g(t_k)) \ge \delta_*/2.
\ee  We can then just apply the theorem to conclude the proof.
 \qed

\medskip

Finally, we give a
\medskip

\noindent {\it Proof of Proposition \ref{prminiSLS}.}

If a minimizer in $W^{1, 2}(M, g(t_0))$ exists for a fixed $t_0$, then argue as in the proof of Theorem \ref{thmrig} (b), we know $(M, g(t))$ is a gradient shrinking soliton. i.e. we can solve the conjugate heat equation with the square of the minimizer as the final value to reach a required integral identity. The proof here is actually easier since $\lam(M, g(t))$ is a constant for solitons.

On the other hand, assume $(M, g)$ is a time slice of a gradient shrinking soliton.
 Then  after normalizing the constant, we can assume, for a smooth function $f$ on $M$, that
\be
\lab{gsseq2}
Ric+ Hess f
-\frac{1}{2} g =0.
\ee
According to  Section 4 in
\cite{CN:1} by Carrillo and Ni, the following equation holds
\be
\lab{cani}
2 \Delta f - |\nabla f |^2 + R + f -n + \mu_s=0
\ee where, for $g^\tau=(1-t) g$, $\tau=1-t$, $t \le 0$ and Perelman's $W$ entropy,
\be
\lab{CNmus}
-\mu_s=\inf_{\tau >0} \mu (g^\tau, \tau) = \inf_{\tau>0}  \inf \{W(g^\tau, u, \tau) \, \, | u \in C^\infty_0(M), \Vert u \Vert_{L^1(M)} =1 \}=\mu(g, 1).
\ee We observe that $-\mu_s$ is just the infimum of the sharp log Sobolev functional in Definition \ref{deflv} i.e.
\be
\lab{mus=lam}
-\mu_s =\lam = \lam(M, g)= \inf \{ L(M, v, g) \, | \, v \in C^\infty_0(M), \, \Vert v \Vert_{L^2(M)} =1 \}.
\ee Here goes the proof. By definition of $\lam \equiv \lam(M, g)$, for any $\e>0$, there exists $v \in C^\infty_0(M)$ with $\Vert v \Vert_{L^2(M)} =1 $ such that
\[
\lam \le L(M, v, g) \le \lam + \e.
\]It is well known by now that, for each fixed $u=v^2$,
\[
\inf_{\tau >0} W(g^\tau, u, \tau) = L(M, v, g).
\]So there exists a $\tau>0$ such that $W(g^\tau, u, \tau) \le L(M, v, g) +\e$. Therefore
\[
\lam \le W(g^\tau, u, \tau) \le \lam + 2 \e.
\]This and \eqref{CNmus} confirms observation \eqref{mus=lam} since $\e>0$ is arbitrary.
Therefore \eqref{cani} becomes
\be
\lab{cani2}
2 \Delta f - |\nabla f |^2 + R + f -n - \lam=0.
\ee Writing
\be
\lab{f2lnw}
f= -2 \ln w- 2 \ln (4 \pi)^{n/4},
\ee then $w$ satisfies
\be
\lab{cani3}
4 \Delta w -R w  + 2 w \ln w  + [\lam  + n + \frac{n}{2} \ln (4 \pi)] w=0.
\ee
It is well known that $f=f(x)$ is comparable to $d^2(0, x)$ when $x$ is large c.f. \cite{CZu:1}. Therefore, from  \eqref{cani3}, $w$ is in $W^{1, 2}(M)$ and it is a minimizer of the sharp log Sobolev functional. This concludes the proof of the proposition.
\qed

\medskip

\medskip
{\bf Acknowledgment.}  We wish to thank Professors Huai-Dong Cao, Pak-Yeung Chan, Ben Chow, Yi Lai and Xialong Li for helpful conversations. We are also grateful to the support of the Simons Foundation through grant No.
710364.

\bigskip

\noindent e-mail:  qizhang@math.ucr.edu

\enddocument